\documentclass[12pt]{article}
\usepackage{epsfig}
\usepackage{amssymb}
\title{Trisecting the $9$-vertex complex projective plane}
\author{Richard Evan Schwartz \thanks{Supported by N.S.F. grant D.M.S.-2102803.}}
\newtheorem{theorem}{Theorem}[section]

\newtheorem{lemma}[theorem]{Lemma}

\def\startproof{{\bf {\medskip}{\noindent}Proof: }}

\def\endproof{$\spadesuit$  \newline}

\def\C{\mbox{\boldmath{$C$}}}%
\def\H{\mbox{\boldmath{$H$}}}%
\def\O{\mbox{\boldmath{$O$}}}%
\def\P{\mbox{\boldmath{$P$}}}%
\def\R{\mbox{\boldmath{$R$}}}%

\begin{document}

\maketitle

\begin{abstract}
  In this paper we will give a short and
  direct proof that
    Wolfgang  K\"{u}hnel's
  $9$-vertex simplicial complex $\C\P^2_9$
  is homeomorphic to $\C\P^2$,
  the complex
  projective plane. The idea
  of our proof is to recall the trisection of
  $\C\P^2$ into $3$ bi-disks and then to
  see this trisection inside a symmetry-breaking
  subdivision of $\C\P^2_9$. After giving
  the proof we will elaborate on the construction
  and sketch an explicit homeomorphism.
  \end{abstract}

  \section{Introduction}

  A $k$-{\it simplex\/} is a $k$-dimensional convex polytope with
  $k+1$ vertices.  For $k=0,1,2,3$ respectively, a $k$-simplex is usually called a
      {\it vertex\/}, {\it edge\/}, {\it triangle\/},  {\it
      tetrahedron\/}.  When $k$ is not important, a $k$-simplex is
    just called a {\it simplex\/}.
    
  A  {\it simplicial complex\/} is
  a finite collection $\cal C$ of simplices, all in an ambient
  Euclidean space, such that
    \begin{itemize}
    \item If $S \in \cal C$ and $S'$ is a sub-simplex of $S$ then $S'
      \in \cal C$.
     \item   If $S,T \in \cal C$ then $S \cap T$ is
       either empty or in $\cal C$.
     \end{itemize}
Informally, the simplices in a simplicial
complex fit together cleanly, without crashing through each other.
     The {\it support\/} $|{\cal C\/}|$ of $\cal C$ is the union
     of all     the  simplices in $\cal C$.   Often we blur the
     distinction between $\cal C$ and $|{\cal C\/}|$ and think of
     a simplicial complex as a union of simplices.
           
  A simplicial complex may be described with no
   mention of the ambient space containing it,
  but there is always the understanding that
    in principle one can find an isomorphic complex in some
    Euclidean space.    To give a pertinent example, let
    $\R\P^2_6$ be the quotient of the
    regular icosahedron by  the antipodal map.  This simplicial complex has
    $6$ vertices, $15$ edges, and $10$ faces.
    One can reconstruct $\R\P^2_6$ in $\R^5$ by fixing some
    $5$-simplex $\Sigma \subset \R^5$, the
    convex hull of vertices $v_1,...,v_6$, then
    mapping
    vertex $k$ of $\R\P^2_6$ to $v_k$ and extending linearly.

\begin{center}
\resizebox{!}{2.2in}{\includegraphics{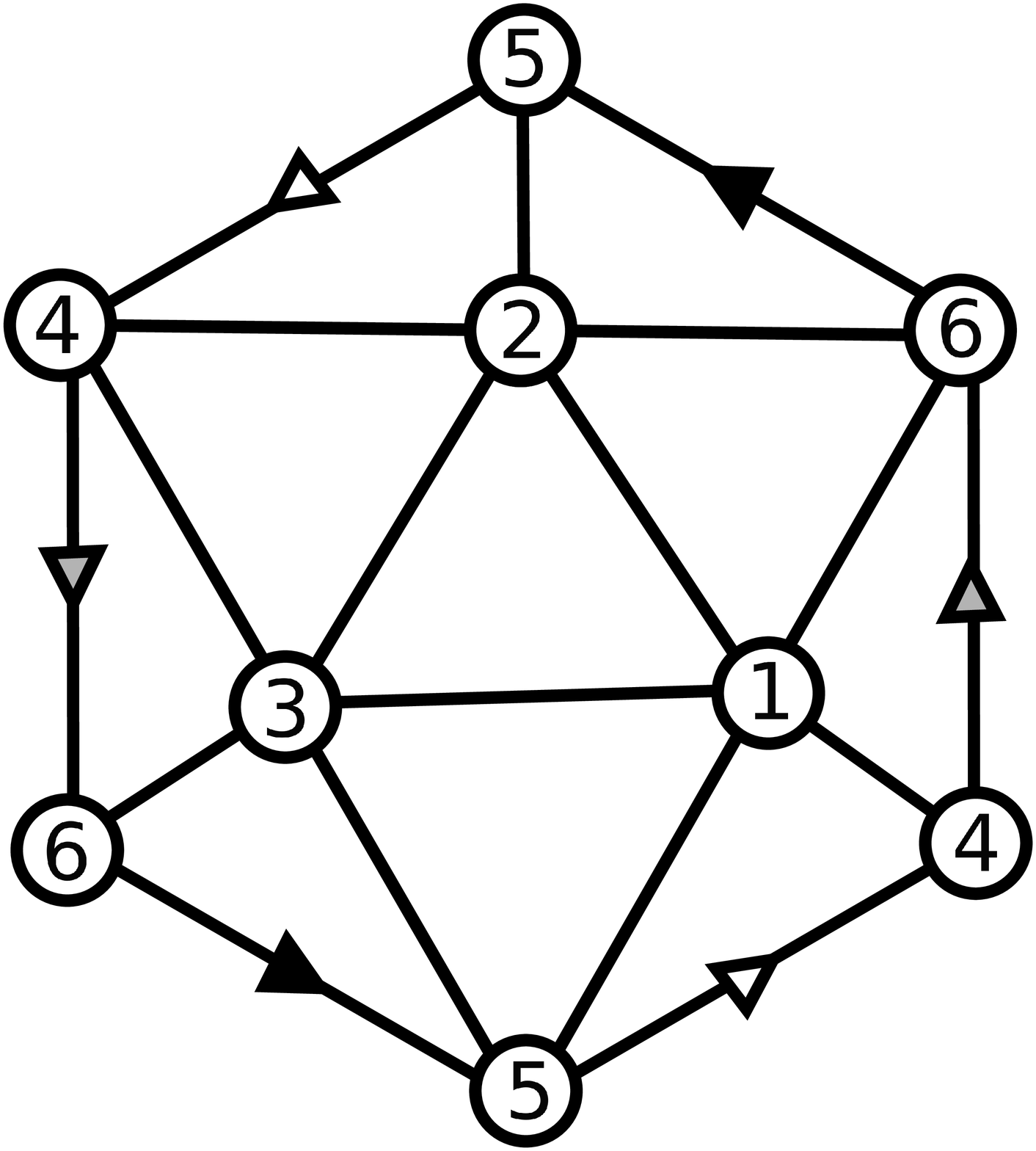}}
\newline
Figure 1: $\R\P^2_6$, the $6$-vertex triangulation of $\R\P^2$.
\end{center}
   
Figure 1 shows another incarnation of
$\R\P^2_6$.  In this picture, the outer edges of the hexagon are
supposed
to be identified according to the labels.
  The complex $\R\P^2_6$ is called a $6$-{\it vertex  triangulation\/}
  of the real projective plane $\R\P^2$ because its support
    is homeomorphic to $\R\P^2$.  This triangulation has
  the fewest number of vertices amongst triangulations of $\R\P^2$, so it is called a
  {\it minimal triangulation\/} of $\R\P^2$.  It is in fact the unique
  minimal triangulation of $\R\P^2$. (Smaller examples like
 the quotient of the regular octahedron by the
  antipodal map fail to be simplicial complexes.)
  
  Here are some other examples related to minimal triangulations.
    \begin{itemize}
\item   The boundary of a tetrahedron is the unique $4$-vertex
  minimal triangulation of the $2$-sphere.  More generally, the
  boundary of a $(k+1)$ simplex is the unique minimal
  triangulation of the $k$-sphere.
\item  If you identify the opposite sides of the big hexagon in
  Figure 4 below, you get the unique minimal
   triangulation $T_7^2$ of the $2$-torus. $T_7^2$ has
$14$ triangles, $21$ edges, and $7$ vertices.
\item   In $1980$,
  W.  K\"{u}hnel discovered $\C\P^2_9$,
  the unique $9$-vertex minimal triangulation of the
  complex projective plane $\C\P^2$.  This triangulation
  has $36$ $4$-simplices and a symmetry group
  of order $54$.
  \item In $1992$, U. Brehm and
  W. K\"{u}hnel [{\bf BK\/}]
  defined $\H\P_{15}^2$ (and two variants), a $15$-vertex
  simplicial complex with $490$ $8$-simplices.  In $2019$,
  D. Gorodkov [{\bf G\/}] proved that $\H\P_{15}$ and the variants
  are PL homeomorphic
  to the quaternionic    projective plane $\H\P^2$.
  \item So far it an open question as to whether there is a
    $27$-vertex
    triangulation of $\O\P^2$, the octonionic (a.k.a. Cayley) projective plane.
  \item The minimal triangulations of $\R\P^3$ and $\R\P^4$
    respectively have $11$ and $16$ vertices. See [{\bf D\/}],
  \item In 2021,
    K. Adiprasito, S. Avvakumov, R. Karasev [{\bf AAK\/}] proved that
   real projective space  can be triangulated
  using a sub-exponential number of simplices.
\end{itemize}  
  The survey article by B. Datta [{\bf D\/}] has a
  wealth of information about minimal triangulations
  up to the year $2007$ and a large number of references.

  The subject of this paper is $\C\P_2^9$.
  In [{\bf KB\/}], 
  K\"{u}hnel and T. Banchoff establish many interesting
  properties of $\C\P_9^2$ and give a rather intricate
  proof that $\C\P_9^2$ really is homeomorphic to
  $\C\P^2$.   Since [{\bf KB\/}], there has been a lot of
  work done trying to understand $\C\P_9^2$ from
  various points of view.  In particular, there are a number of proofs
  that $\C\P_9^2 \cong \C\P^2$, and also a number of proofs that
  $\C\P_9^2$ is the only minimal triangulation of
  $\C\P^2$.  See the article by
  B. Morin and M. Yoshida
  [{\bf MY\/}] for a survey of these proofs.
See also the paper by
B. Bagchi and B. Datta [{\bf BD\/}].

The purpose of this paper is to give a new and very
nice proof that $\C\P_9^2 \cong \C\P^2$.  The basic idea of the proof
here is to recall
the trisection of $\C\P^2$ into $3$ bi-disks, and then to
see this trisection inside a symmetry-breaking subdivision of
$\C\P_9^2$.  The construction is perfectly compatible with an
easier version that works for $\R\P^2_6$, so I will explain that as
well.

The picture developed here is related to the
  $10$-vertex triangulation $\C\P^2_{10}$ of $\C\P^2$ that
  in [{\bf BK\/}] is constructed by
  building outward from $T_7^2$.  Indeed Denis
  Gorodkov, in a private communication, explained to me how
  one can find a ``path'' from
  $\C\P^2_9$ to $\C\P^2_{10}$ using the subdivision idea and then something
  akin to bi-stellar flips. (I'll let Denis tell this story elsewhere
  if he wants to, but see the end of \S 2 for a hint.)
    
  My proof also has a close
  kinship with the ``red-white-blue discussion'' in
    \S 1.3 of the M.P.I.M. preprint by
   Morin and Yoshida that is the precursor to
  [{\bf MY\/}] (and has the same title).
  This discussion is, in turn, related to Figure 8
   in [{\bf KB\/}].  Morin and Yoshida describe the red-white-blue
  discussion as a ``topological insight'' but they don't
  really push it forward into a proof. I think that my
  picture is very similar, but clarified by the special subdivision.

  The approach here possibly could shed light on
    Gorodkov's result that
  $\H\P^2_{15} \cong \H\P^2$.  The same subdivision and
  trisection ideas go through for $\H\P^2_{15}$ almost
  {\it verbatim\/}, and I can see computationally that each of the $3$
  sub-complexes
  is shellable and therefore PL homeomorphic to an $8$-ball.
  However, the high dimensional topology involved in analyzing
  $\H\P_{15}^2$ makes a direct topological analysis of the whole complex
  formidable.  For instance, the sub-complex that plays the role of
  $T^2_7$ has $288$ $6$-simplices.  A key step in extending the proof
  here to $\H\P^2_{15}$ would be showing that this $288$-monster is homeomorphic to
  $(S^3 \times S^3 \times S^3)/S^3$ in a $3$-fold symmetric way.

  Here is an outline of the paper.
  \begin{itemize}
\item   In \S 2 I will give the analogous version of my proof for
  $\R\P^2_6$.  This case is quite concrete and one
  can see the whole idea at a glance.
  \item In \S 3 I will recall the trisection of $\C\P^2$ and discuss
  a few key properties of the {\it central torus\/} in this
  decomposition.
  \item In \S 4 I will describe $\C\P_9^2$ and
  then explain my symmetry-breaking subdivision.
  The construction parallels the real case.
   \item In \S 5 I will
    find the trisection inside the subdivision and
    construct a homeomorphism
    $h: \C\P_9^2 \to \C\P^2$ which
    respects the trisections.
       \item In \S 6 I will
    explain how one can see the real case of
    the construction inside the complex case.
    This analysis leads to a refinement of $h$ and
    gives the full power of our main result,
    Theorem \ref{main}.

  \item In \S 7 I will sketch how to make $h$ completely
    explicit.
    \end{itemize}
    
  I thank Tom Banchoff,  Kenny Blakey, Thomas Goodwillie, Denis Gorodkov,
  Joe Hlavinka, Wolfgang K\"{u}hnel,
  Tyler Lane,
  Dennis Sullivan, and Oleg Viro for helpful discussions.
  (Many of these discussions were about issues related
  to $\H\P_{15}^2$.)
  I also thank the anonymous referee for a number
  of helpful comments, especially those pertaining to
  the real case of the construction.
  These comments from the referee inspired \S 6-7.
  
  \section{The Real Case}

 $\R\P^2$ is the
  space of scale equivalence classes of nonzero 
vectors in $\R^3$.
We denote the equivalence class
of $(x_1,x_2,x_3) \in \R^3$
by $[x_1:x_2:x_3] \in \R\P^2$.

We have the {\it trisection\/}
$\R\P^2=\beta_1 \cup \beta_2 \cup \beta_3$,
where $\beta_j$ is the
set where
$\max(|x_1|,|x_2|,|x_3|)=|x_j|$.
Points in $\beta_1$ may be written
uniquely in the form $[1:x_2:x_3]$, with
$|x_2|, |x_3| \leq 1$.  Thus
$\beta_1$ is a square. So are $\beta_2$ and $\beta_3$.
Each intersection $\beta_i \cap \beta_j$ is a pair of
opposite edges, and the triple intersection
is a union of the $4$ points $[\pm 1:\pm 1:\pm 1]$.
If we interpret $\R\P^2$ as the
quotient of a cube by the antipodal map,
then the $3$ quotient faces are $\beta_1,\beta_2,\beta_3$.

The trisection has $3$-fold symmetry.  The map
$\Sigma: (x_1,x_2,x_3) \to (x_2,x_3,x_1)$
permutes the sets $\beta_1,\beta_2,\beta_3$.
In terms of the cube, $\Sigma$
rotates around the appropriate long diagonal.
$\R\P_6^2$ has a very similar
$3$-fold symmetry: The permutation $S=(123)(456)$ acts as
a rotational symmetry of $\R\P_6^2$.

We add a new vertex $[123]$ at the center of the triangle $(1,2,3)$,
and also new vertices $[12], [13], [23]$ at the centers of the
corresponding edges. 
\begin{center}
\resizebox{!}{1.7in}{\includegraphics{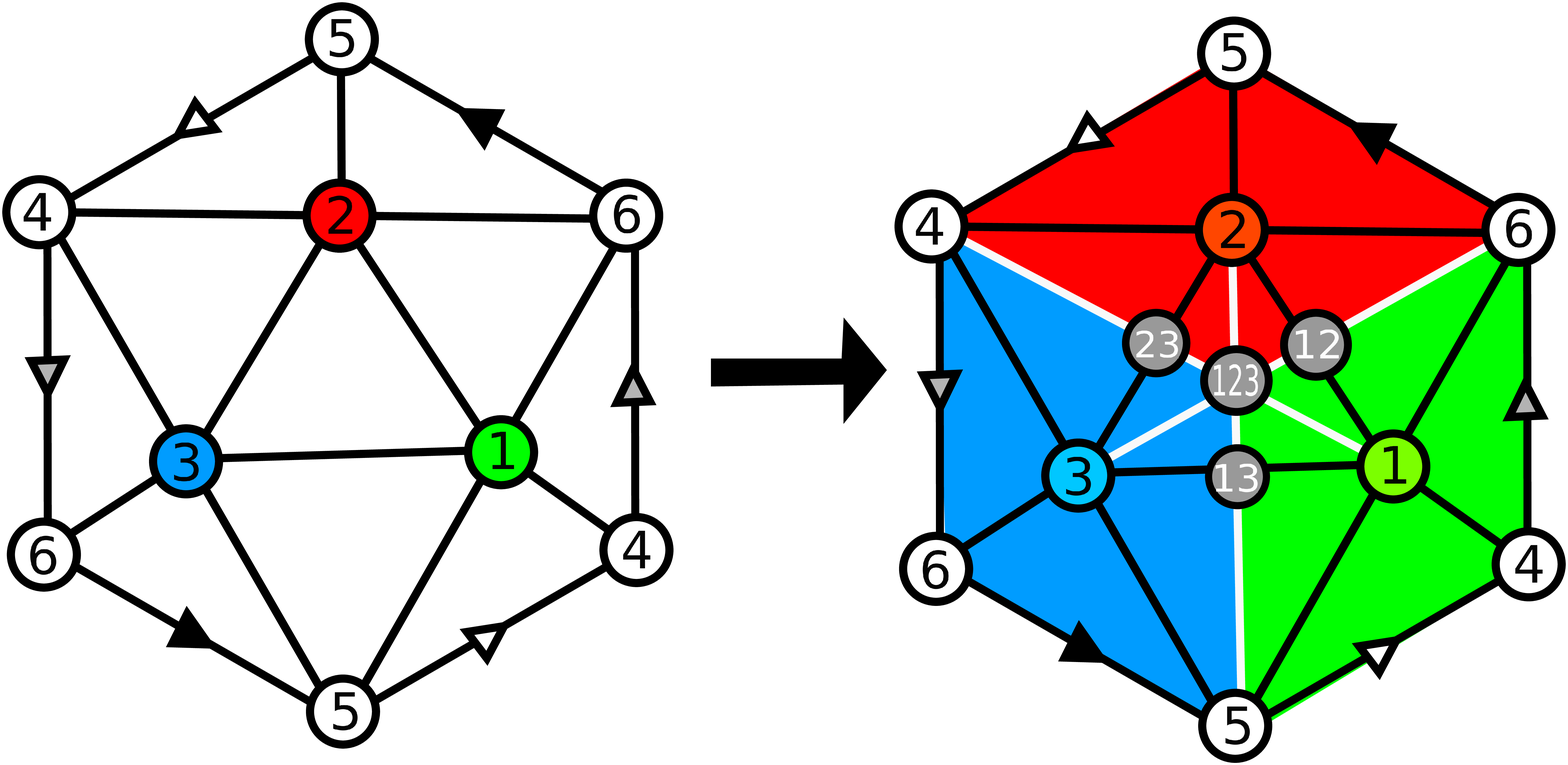}}
\newline
Figure 2: A subdivision of $\R\P^2_6$ into $18$ triangles.
\end{center}
Using the new vertices, we divide the central triangle of
$\R\P^2_6$ into $6$ triangles and we subdivide
each of the adjacent triangles in half.
The subdivision has $18= 3 \times 6$ triangles, with each having
exactly one vertex from the set $\{1,2,3\}$.
For $j=1,2,3$ we let $B_k$ be the subset of $6$ new
triangles having $k$ for a vertex.  The sets
$B_1,B_2,B_3$ are colored green, red, blue in
Figure 2.

This is now the trisection,
and there is a clear homeomorphism from
this subdivided complex to $\R\P^2$ which
maps $B_j$ to $\beta_j$ and conjugates $S$ to $\Sigma$.

Incidentally, a related approach would be to add only $[123]$ and
then to replace the sides $(1,2)$, $(2,3)$, $(3,1)$ with the sides
$([123],6)$, $([123],5)$ and $([123],4)$.  Gorodkov's ``path'' from
$\C\P^2_9$ to $\C\P_2^{10}$ is a more elaborate complex-number analogue of this.

  \section{The Smooth Trisection: Complex Case}

  The complex projective plane $\C\P^2$ is defined
  just as $\R\P^2$ but with respect to the field $\C$
  of complex numbers.
  We denote points in $\C\P^2$
  by $[z_1:z_4:z_7]$. The variable names
will line up with the notation for $\C\P_9^2$.
We have the trisection
$\C\P^2=\beta_1 \cup \beta_4 \cup \beta_7$,
where $\beta_j$ is defined just as in the real case,
using the complex norm in place of the absolute value.
This time, $\beta_j$ is the product of $2$ unit disks.
The bi-disks $\beta_1,\beta_4,\beta_7$
have disjoint interiors and are permuted by the same
map $\Sigma$ as defined in the real case.

The boundary $\partial \beta_1$ is a $3$-sphere, and
it decomposes into the solid tori
$\beta_{14}$ and $\beta_{17}$.  Here
$\beta_{ij}=\beta_i \cap \beta_j$.
To see that $\beta_{14}$ is a solid torus,
note that $\beta_{14}$ consists
of points of the form $[1:u:z]$ with
$|z| \leq |u|=1$ and is therefore the
product of the unit disk and the unit circle.
The {\it central torus\/}
$\beta_{147}=\beta_{14} \cap \beta_{17}=\beta_1 \cap \beta_4 \cap \beta_7$ consists of points where
$|z_1|=|z_4|=|z_7|$.   We discuss $\beta_{147}$ in more detail, with a
view towards seeing it inside $\C\P_9^2$.
\newline
\newline
{\bf Hexagonal Structure:\/}  Let
$\R^3_0 \subset \R^3$ denote the plane of points whose coordinates sum to $0$.
Let $H= \R^3_0 \cap [-1,1]^3$. The vertices of this regular hexagon
are the permutations of $(1,-1,0)$.
 Let $\overline H$ be the
flat torus obtained by identifying the opposite sides
of $H$ by translations. The translation vectors are
the cyclic permutations of $\pm (1,1,-2)$.
 The map
 $$(x_1,x_4,x_7) \to [x^*_1: x^*_4:x^*_7], \hskip 30 pt
x^*=e^{\frac{2 \pi i x}{3}}$$
induces a homeomorphism $\overline H \to \beta_{147}$.
The main point behind this fact
is that  $[1^*:1^*:(-2)^*]=[0^*:0^*:0^*]$, etc.
We equip $\beta_{147}$ with the
metric which makes $\overline H \to \beta_{147}$ an isometry.
\newline
\newline
{\bf Symmetries:\/}
The $3$ fixed points of $\Sigma$ lie in
$\beta_{147}$ and correspond to the points
on $\overline H$ represented by the center and vertices of $H$.
The fixed point set of coordinatewise complex conjugation, which
we call $\cal T$, is $\R\P^2$.  Note that
$\R\P^2 \cap \beta_{147}=\{[\pm 1:\pm 1:\pm 1]\}$.
These points correspond to the center of $H$ and to the
centers of the edges of $H$.
\newline
\newline
{\bf A Contractible Loop:\/}
The line in $\R_0^3$ where $x_1=x_4$ bisects $H$ and contains
the midpoints of a pair of opposite sides.   This line gives
rise to a geodesic loop in $\overline H$.  See the loop
$a_{14}$ in Figure 4 below.  The corresponding loop $\alpha_{14}
\subset \beta_{147}$ is given by $\{[1:1:u]|\ |u|=1\}$.
The loop $\alpha_{14}$ is
contractible in $\beta_{14}$:  It bounds the disk in $\beta_{14}$
consisting of points $[1:1:z]$ with $|z|  \leq 1$.

  \section{The Complex and its Subdivision}

The vertices of $\C\P_9^2$ are labeled $1,...,9$.  Here are 
$16$ of the $36$ $4$-simplices of $\C\P_9^2$ listed on p. 15 of
[{\bf KB\/}].
\begin{itemize}
\item $15289$\ $12389$\ $13689$\ $45289$\ $42389$\ $43689$
\item $14256$\ $14356$\ $14259$\ $14368$
\item $14726$\ $14768$ ($14783$\ $14735$\ $14759$\ $14792$)
\end{itemize}
Comparing our list to [{\bf KB\/}],
we have sometimes permuted the vertices
so as to highlight the indices $1,4,7$.
The other $24$
$4$-simplices are
orbits of the first $12$ under the action of the
{\it fundamental permutation\/}:
$$S=(147)(258)(369).$$
For instance, $14726$ has orbit $14726 \to 14759 \to 14783$.
The four simplices in parentheses
are listed for the sake of making our tetrahedron list below more transparent.
In [{\bf KB\/}] the authors exhibit a symmetry group of
order $54$ acting on $\C\P_9^2$.  For us, one other special element of this group is the
symmetry  $T=(23)(56)(89)$.
\newline

Let $[ij]$ be the midpoint
of the edge $i \leftrightarrow j$.
Let $[ijk]$ be the center of the triangle $ijk$.
Let the {\it rank\/} of a simplex  be the
number of vertices which belong to the set
$\{1,4,7\}$.  Our list above goes by rank.
Parallel to the real case,
we divide each rank $k$ simplex
into $k!$ smaller simplices, as follows:
The rank $1$ simplices are
untouched.
The rank $2$ simplex $14abc$ divides into
$$1[14]abc \hskip 15 pt 4[14]abc,$$
and likewise with the indices $1,4,7$ permuted.
The rank $3$ simplex $147ab$ divides into
{\footnotesize
$$
1[14][147]ab\hskip 10 pt
1[17][147]ab\hskip 10 pt
4[14][147]ab\hskip 10 pt
4[47][147]ab\hskip 10 pt
7[17][147]ab\hskip 10 pt
7[47][147]ab.
$$
\/}
We replace our original $36$ simplices with the subdivided
simplices.   Since there are respectively $18,12,6$ simplices
of rank $1,2,3$ we get a total of $$(1,2,6) \cdot( 18,12,6)=78=3 \times 26$$
{\it new\/} simplices.  (The rank $1$ simplices
count as ``new''.)  Each new simplex has exactly one vertex from the
set $\{1,4,7\}$.

\section{The Combinatorial Trisection}

We have
$\C\P_9^2=B_1 \cup B_4 \cup B_7,$
where
$B_j$ is the union of the
$26$ new simplices having $j \in \{1,4,7\}$ as a vertex.
Each $B_j$ is the cone to vertex $j$  of $\partial B_j$.
Hence $B_i$ and $B_j$ have disjoint interiors for $i \not = j$.
Here is our main result.

\begin{theorem}
  \label{main}
  There is a  homeomorphism $h: \C\P^2_9 \to \C\P^2$ with the following
  properties:
  \begin{enumerate}
\item $h$  maps vertices $1,4,7$ respectively to $[1:0:0]$, $[0:1:0]$, $[0:0:1]$.
\item $h$  maps $B_j$ to $\beta_j$ for $j=1,4,7$.
\item $h$ conjugates $S$ to $\Sigma$.
\item $h$ conjugates  $T$ to $\cal T$.
\end{enumerate}
\end{theorem}
In this section I will construct a non-explicit
homeomorphism $h$ which has the
first
$3$ properties but not necessarily the fourth.
This should satisfy a reader who just wants to see why
$\C\P^2_9 \cong \C\P^2$.  In \S 6, I will give
a more refined version of $h$ which has the fourth property.
In \S 7 I will sketch how to make $h$ explicit.

The first thing we do is list the tetrahedra
in $B_{14}=B_1 \cap B_4.$
We will derive this tetrahedron list from the simplex list above.
The reader might want to check that this actually works, so for
convenience we repeat the simplex list here:
\begin{itemize}
\item $15289$\ $12389$\ $13689$\ $45289$\ $42389$\ $43689$
\item $14256$\ $14356$\ $14259$\ $14368$
\item $14726$\ $14768$ ($14783$\ $14735$\ $14759$\ $14792$)
\end{itemize}

Now for the derivation.
We get $13$ tetrahedra contained in $B_{14}$ by
subdividing the simplices on our list above and omitting $1$ or
$4$.  The tetrahedra are listed in a way that corresponds to the simplices.
\begin{itemize}
\item $5 2 8 9$\hskip 10 pt  $2 3 8 9 $\hskip 10 pt  $3 6 8 9 $
\item $[14] 2 5 6 $\hskip 10 pt  $[14] 3 5 6 $\hskip 10 pt  $[14] 2 5 9$\hskip 10 pt  $[14] 3 6 8 $
\item {\footnotesize $[14] [147] 2 6$\hskip 10 pt    $[14] [147] 6 8 $\hskip 10 pt    
              $[14] [147] 8 3$\hskip 10 pt    $[14] [147] 3 5$\hskip 10 pt    
              $[14] [147] 5 9$\hskip 10 pt    $[14] [147] 9 2 $\/}
          \end{itemize}
The images of these $13$ tetrahedra under $S^{-1}$ lie in $B_{17}$ and
are totally distinct from the ones above. This accounts for all $26$ tetrahedra in $\partial
B_1$.  Hence, the $13$ above are the complete list of tetrahedra
comprising $B_{14}$, and
moreover $B_{14}$ and $B_{17}$ have disjoint interiors.

\begin{lemma}
  \label{lem}
  $B_{14}$ is a solid torus.
\end{lemma}

\startproof
Write $B_{14}=B_{14}' \cup B_{14}''$, where $B_{14}'$ is
the union of the first $3$ tetrahedra above and
$B_{14}''$ is the union of the last ten.
$B_{14}'$ is a $3$-ball because
it is the join of the path $5236$ with the segment
$89$, and $B_{14}''$ is 
a $3$-ball because it is the
cone to vertex $[14]$ of $\partial B_{14}''$, a $10$-triangle
triangulation of the $2$-sphere.

Figure 3 below shows $\partial B_{14}'$ and $\partial B_{14}''$.
Each one is drawn as the union of $2$ combinatorial hexagons
glued along their boundaries according to the labels.
$B_{14}' \cap B_{14}''$ is the union of the $2$ disjoint grey triangles
$259$ and $368$.  Topologically, we get $B_{14}$ by gluing two $3$-balls
together along a pair of disjoint disks in their boundaries.  The
orientations
of the gluings are such that the result is a solid torus (as opposed
to the so-called solid Klein bottle, a nontrivial disk bundle
over the circle).
\endproof

\begin{center}
\resizebox{!}{2.5in}{\includegraphics{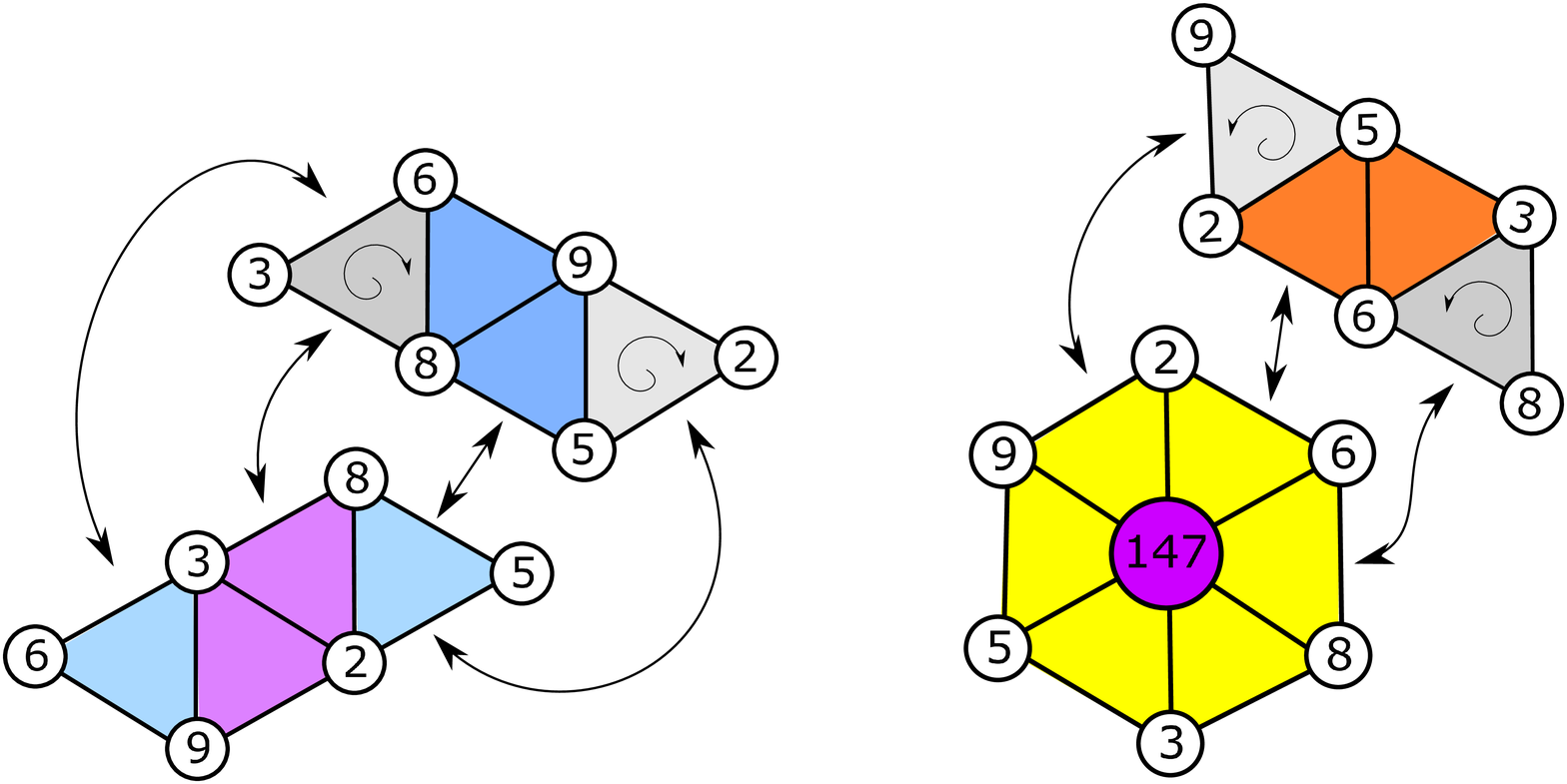}}
\newline
Figure 3: $\partial B_{14}'$ and $\partial B_{14}''$.
Glue the hex boundaries together.
\end{center}

We get the triangulation of $B_{147}=\partial B_{14}$
by gluing the two triangulations from Figure 3 along the
grey triangles.
Figure 4 shows the universal cover of the triangulation.
We get back to $B_{147}$ by gluing the opposite sides of the big hexagon by translations.
 This triangulation of $B_{147}$  is exactly $T_7^2$.
Note that $S$ acts on $B_{147}$ fixing $[147]$, $[258]$, $[369]$,
points which respectively correspond to the center and vertices of
the hexagon, just as in the smooth case.

\begin{center}
\resizebox{!}{3.4in}{\includegraphics{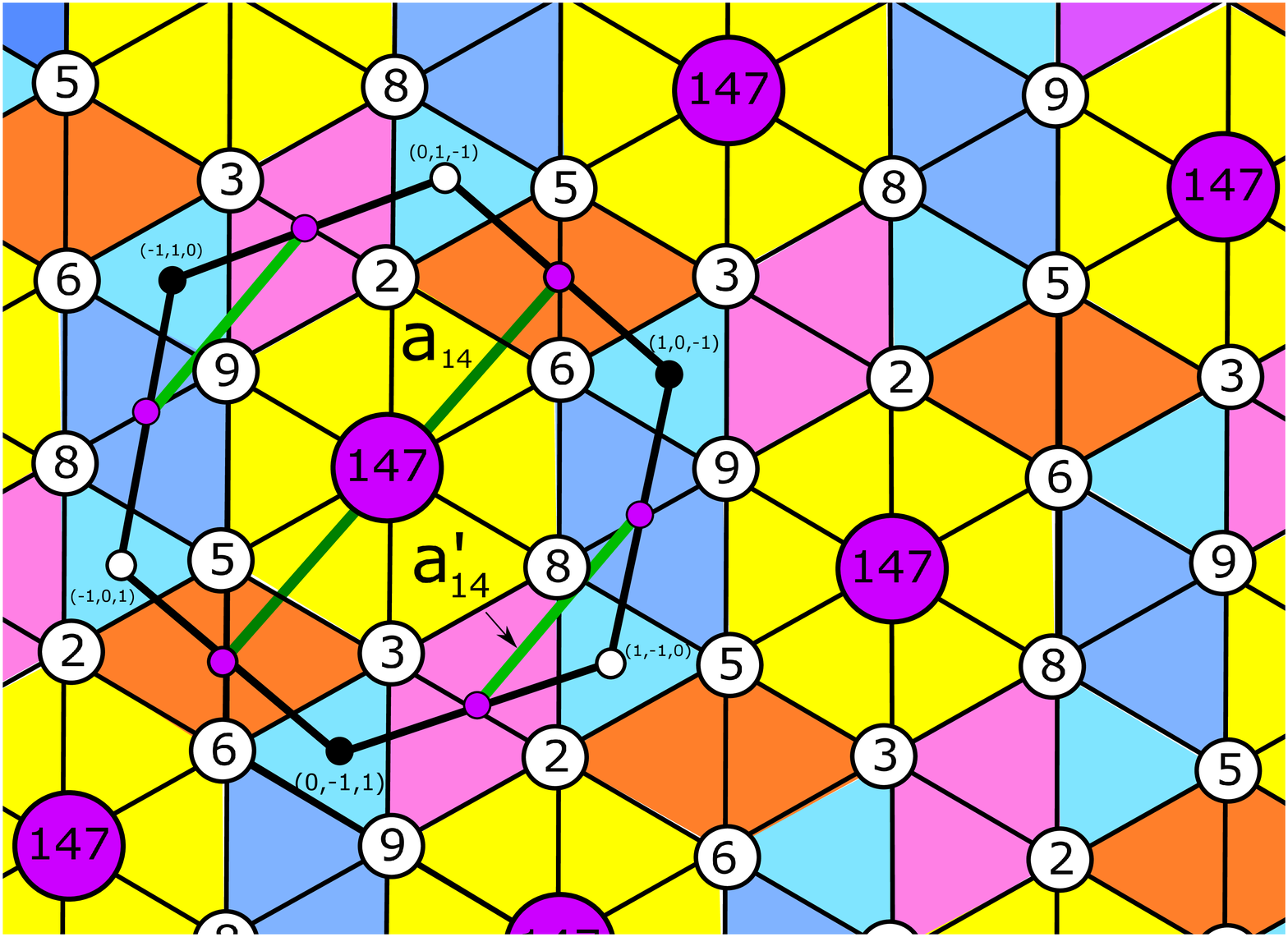}}
\newline
Figure 4: The universal covering of the triangulation of $B_{147}$.
\end{center}

From all this structure we see that (after suitably scaling)
there is an isometry
$h_{147}: B_{147} \to \beta_{147}$ which conjugates $S$ to $\Sigma$
and which maps the green loop $a_{14}$ to $\alpha_{14}$.
The labels of the hexagon vertices, such as $(1,-1,0)$, indicate
precisely how the hexagon here lines up with the one described
in connection with the central torus of $\C\P^2$.
Note that $a_{14}$ is
contractible in $B_{14}$ because $a_{14} \subset B_{14}''$, and
recall that $\alpha_{14}$ is contractible in $\beta_{14}$.
Hence $h_{147}$ extends to 
a homeomorphism $h_{14}: B_{14} \to \beta_{14}$.

Define $h_{17}=\Sigma^{-1} \circ h_{14} \circ S$ and
$h_{47}=\Sigma \circ h_{14} \circ S^{-1}$.  This gives us
homeomorphisms $h_{17}: B_{17} \to \beta_{17}$
and $h_{47}: B_{47} \to \beta_{47}$.  The maps
$h_{ij}$ all agree on $B_{147}$ because $h_{147}$ conjugates
$S$ to $\Sigma$.
The union $$h=h_{14} \cup h_{17} \cup h_{47}: \partial B_1 \cup
\partial B_4 \cup \partial B_7 \to \partial \beta_1 \cup \partial \beta_4 \cup \partial \beta_7$$
is a homeomorphism 
which respects the individual pieces and their intersections.
Since $B_j$ and $\beta_j$ are cones over
$\partial B_j$ and $\partial \beta_j$ we can extend
$h$, by coning, to 
a homeomorphism from
$\C\P_9^2=B_1 \cup B_4 \cup B_7$ to
$\C\P^2=\beta_1 \cup \beta_4 \cup \beta_7$.

\section{The Extra Symmetry}

The fixed set of
$T=(23)(56)(89)$ is a copy of $\R\P_6^2$.  The
$6$ vertices are $1,4,7,[23],[89],[56]$.  If we rename these
vertices $\hat 1,\hat 2,\hat 3,\hat 4, \hat 5,\hat 6$ we get the
same combinatorial pattern as in Figure 1.  Our subdivision
of $\C\P_9^2$ induces the same subdivision as in Figure 2.
In particular, the intersection $\R\P_6^2 \cap B_{14}$ is
a union of $3$ edges which together make $2$ line segments,
namely
$$\hat 4 \hat 5=[23][89] \subset B_{14}', \hskip 30 pt
\hat 6[\hat 1 \hat 2], [\hat 1 \hat 2][\hat 1 \hat 2 \hat 3]=[56][14], [14][147] \subset B_{14}''.$$
The map $h_{147}$ conjugates $T$ to $\cal T$.  Figure 5 below indicates how
$h_{147}$ maps the fixed points of $T$ in $B_{147}$ to the fixed
points of $\cal T$ in $\beta_{147}$.

Now we explain how to choose our homeomorphism $h$ so that it conjugates
$T$ to $\cal T$.
Figure 5 below illustrates the following $4$ disks.
\begin{itemize}
\item Let $D_{14}'$ be the cone to $[[23][89]]$ of the
  loop $a_{14}'$ shown in Figure 4. 
\item Let $D_{14}''$ be the cone to $[14]$ of the loop $a_{14}$ shown in Figure 4.
\item Let $\Delta_{14}'$ be the disk
  $[1:-1:z]$ with $|z| \leq 1$.
  We think of
  $\Delta_{14}'$ as the cone to $[1:-1:0]$ of $\partial \Delta_{14}'$.
\item Let $\Delta_{14}''$ be the disk
$[1:1:z]$ with $|z| \leq 1$.  We think of
$\Delta_{14}''$ as the cone to $[1:1:0]$ of
$\partial \Delta_{14}''$.
\end{itemize}

\begin{center}
\resizebox{!}{1.5in}{\includegraphics{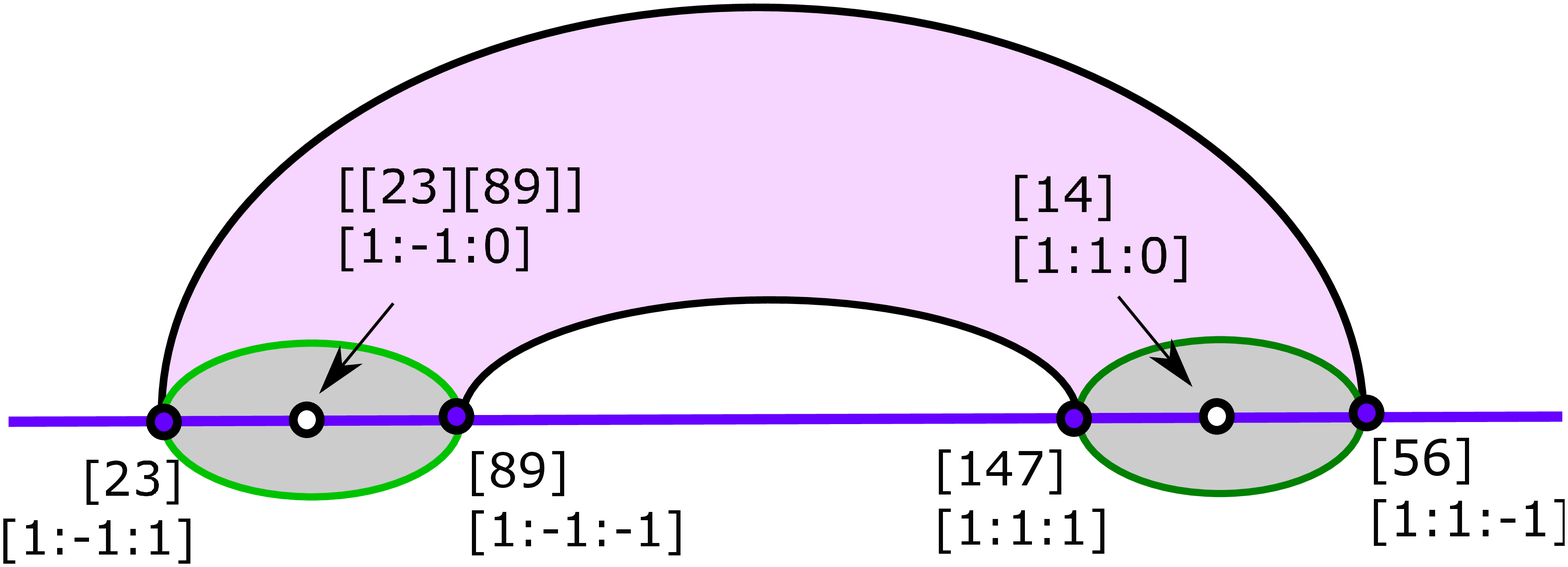}}
\newline
Figure 5: One component of $X_{14}$ or of $\chi_{14}$, depending on the label choice.
\end{center}

Let $Y_{14}$ be the component of $B_{14}-D_{14}'-D_{14}''$ that
contains
the point $[259]$.   Let $\Upsilon_{14}$ be the component of
$\beta_{14}-\Delta'_{14}-\Delta''_{14}$ which contains the point
$[1:1:i]$.   Both $Y_{14}$ and $\Upsilon_{14}$ are solid balls.
$T$ interchanges $Y_{14}$ with (the closure of) $B_{14}-Y_{14}$,
the other component of $B_{14}-D_{14}'-D_{14}''$.
Likewise $\cal T$ interchanges $\Upsilon_{14}$ with (the closure of)
$\beta_{14}-\Upsilon_{14}$.

When we use the top label of each pair, Figure 5 shows $Y_{14}$.
(Think about cutting a pink-frosted  grey donut in half.)
The pink boundary is half of $B_{147}$.
The left grey disk is $D_{14}'$ and the right grey disk is $D_{14}''$.
The map $T$ acts as rotation by $180$ degrees about
the purple line bisecting the grey disks.
The intersection of the purple line
with the grey disks is the part lying in our copy of $\R\P^2_6$.
When we use the bottom labels, Figure 5 shows the same
things for $\Upsilon_{14}$.

By construction $h_{147}(\partial D_{14}')=\partial \Delta_{14}'$ and
$h_{147}(\partial D_{14}'')=\partial \Delta_{14}''$.
We define $h_{14}$ on each of $D_{14}'$ and $D_{14}''$ by coning over the boundaries.
By symmetry this extension conjugates $T$ to $\cal T$ and
is defined in particular on $\partial Y_{14}$.
Our extension maps the (pink and grey) sphere $\partial Y_{14}$ to
the (pink and grey) sphere $\partial \Upsilon_{14}$.
We now extend to a homeomorphism from the ball $Y_{14}$ to
the ball $\Upsilon_{14}$
and use the action of $T$ and $\cal T$ to extend the homeomorphism
to all of $B_{14}$.
Our improved
$h_{14}$ conjugates $T$ to $\cal T$.

The rest of the construction is as above.
  The union map $h_{14} \cup h_{17} \cup h_{47}$, defined
on $\partial B_1 \cup \partial B_4 \cup \partial B_7$,
conjugates $T$ to $\cal T$ because
the pairs $(S,T)$ and $(\Sigma,{\cal T\/})$ commute.
The final coning process respects $T$ and $\cal T$, so the
final extension of $h$ to $B_1 \cup B_4 \cup B_7$
conjugates $T$ to $\cal T$.

Here is where $h$ sends the vertices:
\begin{itemize}
\item $1 \to [1:0:0]$. 
\item $[14]\to [1:1:0]$.
\item $[147] \to [1:1:1]$.
  \item $[56] \to [1:1:-1]$.
\item $2 \to  [1:e^{4 \pi i/7}:e^{12 \pi i/7}].$
\item $[259] \to [1:e^{4 \pi i/7}:0]$.
\end{itemize}
  The remaining images can be readily deduced from the
  action of $S,T,\Sigma, \cal T$.  The last two entries
  require some explanation.
    The coordinates of the point $p_2 \in \R_0^2$ corresponding to
  vertex $2$ are $(-1,5,-4)/7$.   We found this by solving
  the equation
  $2p_2 - \Sigma^2(p_2)=(-1,2,-1)$.
  The choice of where to send $[259]$ is not determined
  by the construction above, but we might as well make it.
  The explicit construction below makes this choice, and so it
  is convenient to list it here.
    
  \section{Making the Homeomorphism Explicit}

  The only non-explicit part of our construction is the extenson of
  the sphere map $h_{147}: \partial Y_{14} \to \partial \Upsilon_{14}$ to
  the ball map $h_{14}: Y_{14} \to \Upsilon_{14}$.  In this section
  we sketch an explicit extension.
  
\begin{center}
\resizebox{!}{2.8in}{\includegraphics{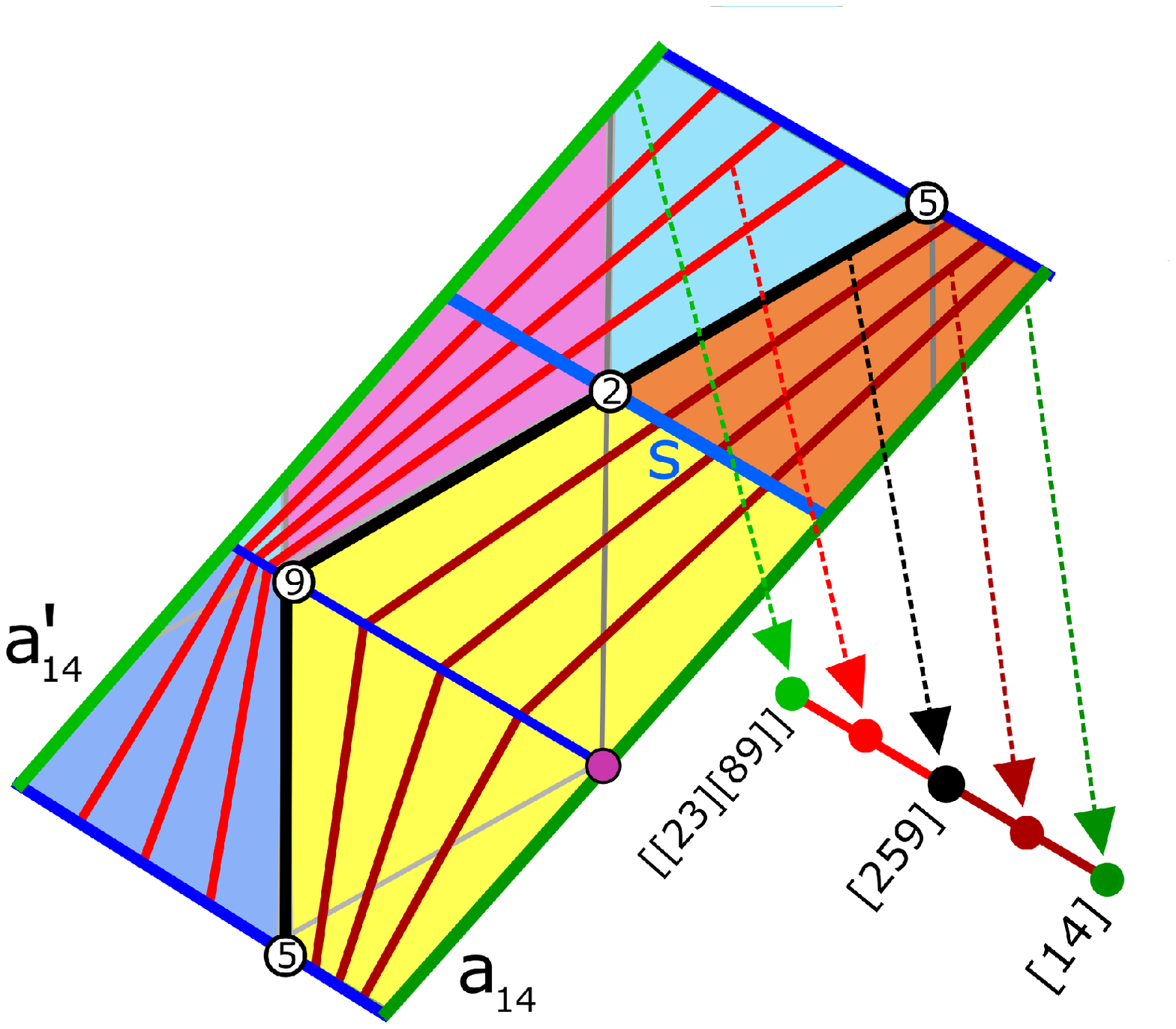}}
\newline
Figure 6: A foliation of $\partial Y_{14}-D_{14}'-D_{14}''$, and the
core of $Y_{14}$.
\end{center}

Gluing the opposite blue sides of the parallelogram
in Figure 6 gives the cylinder $\partial Y_{14}-D_{14}'-D_{14}''$,
drawn pink in Figure 5. The
green loops $a_{14}$ and $a_{14'}$ are on the boundary.
Figure 6 suggests an explicit foliation of $\partial Y_{14}$ by
polygonal loops. Intrinsically these are geodesic bigons in
$B_{147}$.
The intersections with the blue edges of the triangulation move
linearly.

Let the {\it core\/} of $Y_{14}$ be the path with vertices
$[14],[259],[[23][89]]$.
Figure 6 indicates a piecewise linear correspondence between the
loops in the foliation and the points on the core.
  We cone each loop in the foliation to the corresponding point
  on the core.  Now (after some checking of disjointness) we have
  a disk foliation of $Y_{14}$ which interpolates between
 $D_{14}'$ and $D_{14}''$ and respects the partition of $Y_{14}$ into $Y_{14} \cap
 B_{14}'$
 and $Y_{14} \cap B_{14}''$.
  
  We use $h_{147}$ to transfer our
    foliation to $\partial \Upsilon_{14}-\Delta_{14}'-\Delta_{14}''$.
    We change coordinates:
    $$[1:u:z] \to (\arg u,z).$$
    In these coordinates, $\Upsilon_{14} =[0,\pi] \times D^2$.
    Let $s$ be the blue line segment in Figure 6 that
    runs all the way across the cylinder and contains vertex 2.   The
    intersection 
$(t_{\gamma},z_{\gamma})=\gamma \cap s$ is a
    center of symmetry of the loop $\gamma$.
    We cone $\gamma$ to $(t_{\gamma},0)$.
  Now (after some checking of disjointness) we have
  a disk foliation of $\Upsilon_{14}$ which interpolates between
  $\Delta'_{14}$ and $\Delta''_{14}$.  

 We extend $h_{147}$ to $Y_{14}$ by coning, so that it
 maps the one foliation to the other.  This description
 explicitly determines $h_{14}: Y_{14} \to \Upsilon_{14}$.
   Again, the rest of the
  construction is just symmetry and coning, so that final
  map $h$ is explicitly defined.

  \section{References}

  \noindent
  [{\bf AAK\/}], K. Adiprasito, S. Avvakumov, R. Karasev,
  {\it A subexponential
    size  triangulation of $\R\P^n$,\/}
  Combinatorica {\bf 42\/}, 1-8 (2022)
    \newline
    \newline
    \noindent
    [{\bf BD\/}] B. Bagchi and B. Datta, {\it A short proof of the
      uniqueness of  K\"{u}hnel's $9$-vertex complex projective
      plane\/}, Advances in Geometry (2001) pp 157-163.
    \newline
    \newline
    [{\bf BK\/}] T. Banchoff and W  K\"{u}hnel,
    {\it Equilibrium triangulations of the complex projective
      plane\/},
    Geometriae Dedicata {\bf 44\/} (1992) pp 311-333
    \newline
    \newline
    [{\bf BrK\/}] U. Brehm and W  K\"{u}hnel,
    {\it $15$-vertex triangulations of an $8$-manifold\/},
    Math Annalen {\bf 294\/} (1992) pp 167--193
    \newline
    \newline
      [{\bf D\/}] B. Datta, {\it Minimal Triangulations of
        Manifolds\/},
      Journal of the Indian Institute of Science, Vol. 87, No. 4,
      (2007)
      pp. 429-450.  \newline
      (See also
  arXiv:math/0701735v1)
  \newline
  \newline
  [{\bf G\/}] D. Gorodkov, {\it A $15$-Vertex Triangulation of the
    Quaternionic Projective Plane\/}, Discrete Computational Geometry,
 (Sept 2019) pp 348-373
  \newline
  \newline
[{\bf KB\/}]: W. K\"{u}hnel and T. F. Banchoff,
{\it The $9$-Vertex Complex Projective Plane\/},
Math. Intelligencer.  Vol 5, No 3 (1983) pp 11-22
\newline
\newline
[{\bf MY\/}] B. Morin and M. Yoshida, {\it The
      K\"{u}hnel  trangulation of the complex projective
      plane from the point of view of complex crystallography\/},
    Mem, Fac. Sci. Kyushu Univ. Ser. A {\bf 45\/} (1991) pp 55-142
    \noindent

\end{document}